\documentclass[10pt]{article}%
\usepackage{amsmath}
\usepackage{amsfonts}
\usepackage{amssymb}
\usepackage[all]{xypic}
\usepackage{hyperref}
\usepackage{amsmath,amsfonts,amsthm}
\input cyracc.def
\font \tencyr=wncyr10
\newfam\cyrfam
\font\tencyr=wncyr10
\font\sevencyr=wncyr7
\font\fivecyr=wncyr5
\def\cyr{\fam\cyrfam\tencyr\cyracc}

\textfont\cyrfam=\tencyr \scriptfont\cyrfam=\sevencyr
\scriptscriptfont\cyrfam=\fivecyr
\let\kappa\varkappa
\newcommand{\ev}{\text{\cyr \`{E}}}

\newcommand{\N}{\mathbb{N}}

\newcommand{\F}{\mathcal{F}}
\newcommand{\C}{\mathcal{C}}
\newcommand{\E}{\mathcal{E}}
\newcommand{\NN}{{\partial B}}
\newcommand{\LL}{\mathcal{L}}

\newcommand{\Diff}{\mathrm{Diff\,}}
\newcommand{\im}{\mathrm{im\,}}

\newcommand{\df}{\stackrel{\mathrm{def.}}{=}}
\newcommand{\REF}[1]{{\normalfont (\ref{#1})}}
{\theoremstyle{remark}\newtheorem{remark}{Remark}}
{\theoremstyle{remark}}
\newtheorem{theorem}{Theorem}

\newtheorem{corollary}{Corollary}
\newtheorem{definition}{Definition}
\newtheorem{proposition}{Proposition}
\begin{document}
\date{\today}
\title{Domains in Infinite Jets: $\mathcal{C}$--Spectral Sequence}
\author{A.~M. Vinogradov\thanks{Dipartimento di Matematica e Informatica, Universit\`{a} degli Studi di
Salerno, Via Ponte don Melillo, 84084 Fisciano (SA), Italy;
INFN, Gruppo Collegato di Salerno, Italy; email: vinograd@unisa.it.} and G. Moreno\thanks{Dipartimento Matematica e Applicazioni ``Renato Caccioppoli'', Universit\`{a} degli Studi di
Napoli, Via Cintia, 80126 Napoli (NA), Italy; email: giovanni.moreno@dma.unina.it; giovanni.moreno@fastwebnet.it}}
\maketitle

\begin{abstract}
Domains in infinite jets present the simplest class of diffieties with boundary.
In this note some basic elements of geometry of these domains are introduced and
an analogue of the $\mathcal{C}$--spectral sequence (see \cite{CSpecSeqI,CSpecSeqII})
in this context is studied. This, in particular, allows cohomological interpretation
and analysis in the spirit of  \cite{CoAn} of initial data, boundary conditions, etc,
for general partial differential equations and of transversality conditions in calculus
of variations. This kind applications and extensions to arbitrary diffieties will be
considered in subsequent publications.
\end{abstract}

\section*{Introduction}
The diffiety representing a system of PDEs $\E$ is its infinite
prolongation $\E_{\infty}$, while  an initial data problem
associated with $\E$ is represented by a subdiffiety of
$\E_{\infty}$ which is, in a sense, the infinite prolongation of
the original initial data and whose codimension and co-Dimension
are both equal to $1$ (see \cite{CSpecSeqII}). This is one of numerous
situations when the necessity to study a pair of diffieties and,
in particular, a diffiety with boundary arises. Such a pair will
be denoted $(B,\partial B)$ even when $\partial B$ is not the
boundary of $B$. The simplest situation of this kind is associated
with a pair of smooth manifolds $(E,\partial E)$, $\partial E$
being a hypersurface in $E$, if one puts $B=J^\infty(E,n)$ and
\begin{equation}\label{defDeBi}
\partial B\df\{[L]_y^\infty\ |L\subset E, \dim L=n, \ L\textrm{ intersects }
\partial E\textrm{ transversally at }y\}.
\end{equation}
In the case when  $E\stackrel{\pi}{\longrightarrow}M$ is a smooth
fiber bundle over an $n$--dimensional manifold $M$ and $\partial
E=\pi^{-1}(\partial M)$, with $\partial M$ being a hypersurface in
$M$, the above construction applied to graphs of sections of $\pi$ gives $B=J^\infty(\pi)$ and $\partial
B=\pi_\infty^{-1}(\partial M)$. This particular case will be
referred to as \emph{fibered}. In the most general case, $B$ is an
open domain in $J^\infty(E,n)$ and $\partial B$ is
\emph{$\partial$--admissible} in $B$ (see \cite{CSpecSeqII}, 8.4).
\par In this note the \emph{relative $\C$--spectral sequence of the pair
$(B,\partial B)$}, denoted by $(E_r(B,\partial
B),d_{r,\textrm{rel}})$, will be constructed by following the
guidelines of \cite{CSpecSeqII} (see 12.5).
\par In the sequel we follow the notation of
\cite{CoAn}. Namely, $\F$ and $\Lambda$ stand for the filtered
algebra of smooth functions and for the algebra of differential
forms on $B$, respectively, $\C\subset\Lambda$ for the ideal of
Cartan forms, $\overline{\Lambda}$ for the differential algebra of
horizontal forms and $\kappa$ for the Lie algebra of higher
symmetries of $B$. If $B=J^\infty(\pi)$, then the evolutionary
derivation, whose generating function is $\psi$, is denoted by
$\ev_\psi$, $\ell_f$ stands for the universal linearization of
$f\in \F$, and so on. Accordingly, symbols $\F(\NN)$,
$\Lambda(\NN)$, $\C_\NN$, $\overline{\Lambda}(\NN)$,
$\kappa(\NN)$, $\ev_\psi^\NN$, and $\ell_f^\NN$ stand for the
corresponding objects on the diffiety $\NN$. For instance,
$\ell_f^\NN$ is a $\C_\NN$--differential operator. The
$\C$--spectral sequence for $B$ is denoted by $E_r$, while the
$\C_\NN$--spectral sequence for $\NN$  by  $E_r(\NN)$.

A suitable for describing the above situation local chart
$(x_1,\ldots,x_n)$ on $M$ is such that $\partial M=\{x_n=0\}$. It
is extended to a chart on $E$ by introducing some fiber
coordinates $(u^1,\ldots,u^m)$ and then to the standard jet chart
$(x_1,\ldots,x_n,u^1,\ldots,u^m,\ldots,u^k_{\sigma},\ldots)$ on
$J^\infty(\pi)$. $\partial B$ in this chart is given by
$\{x_n=0\}$. Total derivatives $D_i, \;1\leq i\leq n$, on
$J^\infty(\pi)$ corresponding to this chart are tangent to
$\partial B$ if $i<n$. Denote by $\Pi^{(j)}, \;j=1,\ldots,m$, the
projection of $\F(\pi,\pi)$ to its $j$--th component and put
$D_\sigma^{(j)}=D_\sigma\circ \Pi^{(j)}$ with
$\sigma=(\sigma_1,\ldots,\sigma_n)\in\N_0^n$ being a multi--index.
Elements of the standard basis in $\N_0^n$ are denoted by $1_i,
\;i=1,\ldots,n$. By fixing volume forms $dx_1\wedge\cdots\wedge
dx_n$ and $dx_1\wedge\cdots\wedge dx_{n-1}$ on $M$ and $\partial
M$, respectively, we identify $\overline{\Lambda}^n$--valued
operators  with $\F$--valued and
$\overline{\Lambda}^{n-1}(\NN)$--valued operators with
$\F(\NN)$--valued ones, respectively.

\section{The relative $\C$--spectral sequence, general case}

Let $\iota_\NN:\NN\rightarrow B$ be the embedding map. Denote by
$\LL=\Lambda(B,\NN)$ the ideal $\ker \iota_\NN^*$ of vanishing on
$\partial B$ differential forms on $B$. Then the restriction of
the Cartan distribution on $B$ to $\NN$ is given by the ideal
\begin{equation}
\C_\NN=\frac{\C}{\C\cap\mathcal{L}}
\end{equation}
of the quotient algebra $\Lambda/\LL=\Lambda(\NN)$.
\begin{proposition}
The submodule
\begin{equation}
E_0^p(B,\NN)\df \frac{\mathcal{C}^p\cap\mathcal{L}+\mathcal{C}^{p+1}}{\mathcal{C}^{p+1}}
\end{equation}
of $E_0^p$ is, moreover, a sub--complex of $(E_0^p,\ d_0)$.
\end{proposition}
\begin{definition}
The term $E_0^0(B,\NN)$ is called the differential algebra of
\emph{relative with respect to $\partial B$ horizontal
differential forms on $B$} and is denoted by $\overline{\Lambda}
(B,\,\NN)$.
\end{definition}
\begin{theorem}\label{teoE0}
The quotient $\F$--module $\frac{E_0^p}{ E_0^p(B,\NN)}$ is
isomorphic to the $\F(\NN)$--module
$\frac{\C_\NN^p}{\C_\NN^{p+1}}$, i.e., to $E_0^p(\NN)$.
\end{theorem}
This way one gets the short exact sequence of complexes
\begin{equation}\label{eqShExSeqE0}
0\to  E_0^p(B,\NN)\stackrel{i}{\to} E_0^p\stackrel{\alpha}{\to}
E_0^p(\NN)\to 0,
\end{equation}
which leads to the corresponding long exact cohomology sequence
\begin{equation}\label{eq146'}
\xymatrix{
E_1^p(B,\NN) \ar[rr]^{H(i)} && E^p_1\ar[dl]^{H(\alpha)}\\
& E^p_1(\NN)\ar[ul]^{\partial}}
\end{equation}
where
\begin{equation}
E_1(B,\NN) \df H(E_0(B,\NN)).
\end{equation}
\begin{proposition}\label{propOneLine}
If the the one--line theorem (see \cite{CoAn}, 4.3.7) holds for
both $B$ and $\partial B$, then it holds also for the relative
$\C$--spectral sequence of the pair $(B,\partial B)$.
\end{proposition}
In other words, if the term $E_1^{p,q}(B,\NN)$ is nontrivial, then
either $p=0$, or $q=n$.
\par In particular, if $p=0$,  then \REF{eq146'} reads
\begin{equation}\label{eq148}
\xymatrix{
\overline{H}(B,\NN) \ar[rr]^{H(i)} && \overline{H}\ar[dl]^{H(\alpha)}\\
& \overline{H}(\NN)\ar[ul]^{\overline{\partial}}}
\end{equation}
and is called the \emph{long exact sequence of horizontal de Rham
cohomologies} of the pair $(B,\NN)$. Here $\overline{\partial}$ is
the \emph{horizontal coboundary operator}. Under the hypothesis of
Proposition \ref{propOneLine}, sequence \REF{eq146'} for $p>0$
reduces to the short exact sequence
\begin{equation}\label{eqSeqEsCorE1}
0
{\longrightarrow}E_1^{p,n-1}(\NN)\stackrel{\partial}{\longrightarrow}E_1^{p,n}
(B,\NN)\stackrel{H(i)}{\longrightarrow}E_1^{p,n}
{\longrightarrow}0 .
\end{equation}

\section{The relative $\C$--spectral sequence, fibered case}
\subsection{Description of $\partial B$}
Let $\Delta$ be a completely integrable distribution on $M$.
Denote by $M_x$ the leaf of the corresponding to $\Delta$
foliation that passes through $x\in M$.
\begin{definition}
Two (local) sections $s,s'\in\Gamma_{\mathrm{loc}}(\pi)$ are said
to be $k$--th order tangent along $\Delta$ at the point $x$ if
their restrictions $s_{|M_x},s'_{|M_x}$ to   $M_x$ are $k$--th
order tangent at $x$ in the usual sense.
\end{definition}
Denote by $[s]_{\Delta,x}^k$ the represented by $s$ equivalence
class of local sections of $\pi$  that are $k$--th order tangent
each other along $\Delta$ at $x$ and put
\begin{equation}
J^k_\Delta(\pi) \df \{[s]_{\Delta,x}^k\ |\ s\in \Gamma_{\mathrm{loc}}(\pi),\ x\in M\},
\end{equation}
\begin{equation}
\pi^\Delta_{k,l}: J^k_\Delta(\pi) \to J^l_\Delta(\pi),\quad \pi^\Delta_{k,l}\left( [s]_{\Delta,x}^k \right)\df [s]_{\Delta,x}^l ,\quad k\geq l,
\end{equation}
\begin{equation}
\pi^\Delta_{k}: J^k_\Delta(\pi) \to M,\quad \pi^\Delta_{k}\left( [s]_{\Delta,x}^k \right)\df x.
\end{equation}
The inverse limit $ J^\infty_\Delta(\pi)\stackrel{\pi^\Delta_\infty}{\longrightarrow}M$ of $\pi_k^\Delta$, $k\to\infty$, is defined in the usual way.\par
The following results clarify how $ J^\infty_\Delta(\pi) $ is related to $ J^\infty(\pi)$ and the infinite jets bundle of the restriction $\pi_{|M_x}$.
\begin{proposition}\label{propDueMappe}
The maps
\begin{equation}\label{map1}
J^\infty(\pi_{|M_x})\to J^\infty_\Delta(\pi),\quad   [s]_x^\infty\mapsto [\widetilde{s}]_{\Delta,x}^\infty,
\end{equation}
with $\widetilde{s}$ being an extension of $s\in\Gamma_{\mathrm{loc}}(\pi_{|M_x})$ to $\Gamma_{\mathrm{loc}}(\pi)$, and
\begin{equation}\label{map2}
J^\infty(\pi) \to J^\infty_\Delta(\pi),\quad  [s]_x^\infty \mapsto [s]_{\Delta,x}^\infty,
\end{equation}
are injective and surjective, respectively.
\end{proposition}
\begin{remark}
Observe that $(M,\Delta)$ is a diffiety so that it makes sense to
consider $\Delta$--differential operators. If $\pi$ is linear, the
sub--functor $\Delta\Diff_k(\Gamma(\pi),\,\cdot\,)$ of
$\Diff_k(\Gamma(\pi),\,\cdot\,)$ is represented by
$\mathcal{J}^k_\Delta(\pi)=\Gamma(\pi_\Delta^k)$. The  projection
of representative objects
$\mathcal{J}^k(\pi)\mapsto\mathcal{J}^k_\Delta(\pi)$ corresponds
to the natural inclusion of functors
$\Delta\Diff_k(\Gamma(\pi),\,\cdot\,)\subset\Diff_k(\Gamma(\pi),\,\cdot\,)$.
This is the meaning of \REF{map2}. On the other hand, since
$\Delta$--differential operators admit restrictions to the leaves
of $\Delta$, the natural projection of functors
$\Delta\Diff_k(\Gamma(\pi),\,\cdot\,)\to
\Diff_k(\Gamma(\pi_{|M_x}),\,\cdot\,)$ is represented by an
injection
$\mathcal{J}^k(\pi_{|M_x})\subset\mathcal{J}^k_\Delta(\pi)$ of the
corresponding representative objects.  This is the meaning of
\REF{map1}.
\end{remark}
Now, let $\nabla$ be a complementary to $\Delta$  completely
integrable distribution on $M$. Then $D(\Delta)\oplus
D(\nabla)=D(M)$. Assume also that $\partial M$ is a leaf of
$\Delta$.

\begin{proposition}\label{propIdentificazione}
Diffieties $J^\infty_\nabla(\pi^\Delta_\infty)$ and $J^\infty_\Delta(\pi^\nabla_\infty)$
are both identified naturally with $J^\infty(\pi)$.
\end{proposition}
Put $\pi^\nabla_{\infty,\partial
M}\df(\pi^\nabla_\infty)_{|\partial M}$ and, by using Propositions
\ref{propDueMappe} and \ref{propIdentificazione}, define the
embedding
\begin{equation}
J^\infty(\pi^\nabla_{\infty,\partial M})\subset  J^\infty_\Delta(\pi^\nabla_\infty)=
J^\infty(\pi).
\end{equation}
Fiber bundle $\pi^\nabla_{\infty,\partial M}$ will be referred to as the (infinite)
\emph{normal jets bundle} of $\pi$ with respect to the hypersurface $\partial M$.
The following result sheds light on the structure of $\partial B$ in the fibered case.
\begin{theorem}
Let $\Delta$,  $\nabla$  and $\partial M$ be as above. Then
diffieties $\NN$ and $J^\infty(\pi^\nabla_{\infty,\partial M})$
are naturally identified.
\end{theorem}
A section of $\pi_\infty^\nabla$ is described by a vector $\boldsymbol{f}=(\ldots,f_i^k,\ldots)$, $f_i^k\in C^\infty(M)$, and, correspondingly, a section of $\pi^\nabla_{\infty,\partial M}$ is described by a similar vector, with the $f_i^k\in C^\infty(\partial M)$. So, an element $\theta$ of $J^\infty(\pi^\nabla_{\infty,\partial M})$ is represented by the vector $(x_1,\ldots,x_{n-1},\ldots,\frac{\partial^{|\tau|} f_i^k}{\partial x^\tau}(x_1,\ldots,x_{n-1}),\ldots)$, with $k=1,\ldots,m$, $i\in\N_0$, and $\tau\in\N_0^{n-1}$. Put
\begin{equation}
(u_i^k)_\tau(\theta)\df \frac{\partial^{|\tau|} f_i^k}{\partial x^\tau}(x_1,\ldots,x_{n-1}),\quad\theta=[\boldsymbol{f}]_{(x_1,\ldots,x_{n-1})}^\infty.
\end{equation}
The embedding $\iota_\NN$ in the considered case is described by
\begin{equation}
\left\{\begin{array}{rcll}\iota_\NN^*(x_i)&=&x_i, &  i=1,\ldots,n-1, \\\iota_\NN^*(x_n)&=&0, &  \\\iota_\NN^*(u_\sigma^k)&=&(u_{\sigma_n}^k)_{\sigma-\sigma_n1_n}, &  k=1,\ldots,m,\sigma\in\N_0^n.\end{array}\right.
\end{equation}

\subsection{The $E_0$ term of the relative $\C$--spectral sequence}
Assume that $\NN$ is given by the equation $\phi=0$,
$\phi=\pi_{\infty}^*(\varphi)$. Then  the following relations are
easily checked:
\begin{eqnarray}
\Lambda(B,\NN)&=&\phi \Lambda + d\phi\wedge\Lambda,\\
\overline{\Lambda} (B,\,\NN)&=&\phi \overline{\Lambda} +\overline{d}
\phi\wedge \overline{\Lambda},\\
E_0^p(B,\NN)&=&\phi E_0^p+\overline{d}\phi\wedge E_0^p.\label{eq19'}
\end{eqnarray}
\begin{proposition}\label{propE0rel}
The following $\F$--modules isomorphism holds:
\begin{equation}\label{eq77}
E_0^p(B,\NN)\cong\C\Diff_{(p)}^\mathrm{alt.}(\kappa,\,\overline{\Lambda} (B,\,\NN)).
\end{equation}
\end{proposition}
\begin{corollary}\label{corE0}
The quotient $\frac{E_0^p}{ E_0^p(B,\NN)}$ is isomorphic to
$\C_\NN\Diff_{(p)}^\mathrm{alt.}(\kappa(\NN),\,\overline{\Lambda}(\NN))$.
\end{corollary}
It is not difficult to see that the module
$\Gamma(\pi_{\infty,\partial M}^\nabla)$ is locally free. So, an
element $\boldsymbol{\psi}$ of $\F(\NN,\pi_{\infty,\partial
M}^\nabla)$ can be represented by its generating function
$\boldsymbol{\psi}=(\ldots,\psi_i^k,\ldots)$,
$\psi_i^k\in\F(\NN)$, with $k=1,\ldots,m$ and $i\in\N_0$. The
higher symmetry of $\NN$ corresponding to $\boldsymbol{\psi}$ is
represented by the evolutionary derivation
\begin{equation}
\ev_{\boldsymbol{\psi}}^\NN\df\sum D_\tau(\psi_i^k)\frac{\partial}
{\partial(u_i^k)_\tau},\quad \tau\in\N_0^{n-1}, i\in\N_0, k=1,\ldots,m.
\end{equation}
Denote by $\Pi^{(k,i)}$ the projection of the free module
$\F(\NN,\pi_{\infty,\partial M}^\nabla)$ onto its $(k,i)$--th component and
put $D_\tau^{(k,i)}=D_\tau\circ \Pi^{(k,i)}$. The above formula now reads
\begin{equation}
\ev_{\boldsymbol{\psi}}^\NN\df\sum D^{(k,i)}_\tau(\boldsymbol{\psi})
\frac{\partial}{\partial(u_i^k)_\tau},\quad \tau\in\N_0^{n-1}, i\in\N_0,
k=1,\ldots,m,
\end{equation}
and allows to introduce the universal linearization operator
\begin{equation}
\ell_g^\NN\df \sum \frac{\partial g}{\partial(u_i^k)_\tau}D_\tau^{(k,i)},
\quad \tau\in\N_0^{n-1}, i\in\N_0, k=1,\ldots,m,\quad g\in\F(\NN).
\end{equation}
\begin{proposition}
The isomorphism of Corollary \ref{corE0} between
$\frac{E_0^{1,0}}{E_0^{1,0}(B,\NN)}$ and\linebreak 
$\C_\NN\Diff(\kappa(\NN),\,\F(\NN))$ sends the equivalence class
$[{U}_1(f)]$ to the $\C_\NN$--differen\-tial operator
$\ell_{f_{|\NN}}^\NN$.
\end{proposition}
In view of Proposition \ref{propE0rel} and Corollary \ref{teoE0}
sequence \REF{eqShExSeqE0} is identified with
\begin{equation}\label{eqShExSeqE0bis}
0\to  \C\Diff_{(p)}^\mathrm{alt.}(\kappa,\,\overline{\Lambda} (B,\,\NN))
\stackrel{i}{\to}  \C\Diff_{(p)}^\mathrm{alt.}(\kappa,\,\overline{\Lambda})
\stackrel{\boldsymbol{\alpha}}{\to} \C_\NN\Diff_{(p)}^\mathrm{alt.}(
\kappa(\NN),\,\overline{\Lambda}(\NN))\to 0.
\end{equation}
\begin{proposition}\label{prop12}
The projection $\boldsymbol{\alpha}$ in \REF{eqShExSeqE0bis} corresponding to
$\alpha$ in \REF{eqShExSeqE0}  is given by
\begin{equation}
\boldsymbol{\alpha}(\ell_{f^1}\wedge\cdots\wedge\ell_{f^p}\otimes\overline{\omega})=
\ell_{f^1_{|\NN}}^\NN\wedge\cdots\wedge\ell_{f^p_{|\NN}}^\NN\otimes\overline{\iota_\NN^*(\omega)}
\end{equation}
with $f^1,\ldots,f^p\in\F$ and $\overline{\omega}\in\overline{\Lambda}$.
\end{proposition}
\subsection{The $E_1$ term of the relative $\C$--spectral sequence}
Proposition \ref{propOneLine} holds in the fibered case, so sequence \REF{eqSeqEsCorE1} reads
\begin{equation}\label{eqSeqEsCorE1bis}
0\to \widehat{\kappa}(\NN)\stackrel{\partial}{\longrightarrow}\widehat{\kappa}(B,\NN)\stackrel{H(i)}{\longrightarrow}\widehat{\kappa}\to 0.
\end{equation}
The \emph{relative adjoint} to $\kappa$ module  is defined as
\begin{equation}
\widehat{\kappa}(B,\NN)\df E_1^{1,n}(B,\NN),
\end{equation}
and the module product in it is defined by composing
$\C$--differential operators with scalars from the right. With
respect to this module structure both $\partial$ and $H(i)$ are
$\F$--linear maps.
\par
Since $\overline{\Lambda}^n(B,\NN)=\overline{\Lambda}^n$,
\begin{equation}
\widehat{\kappa}(B,\NN)=\frac{\C\Diff(\kappa,\overline{\Lambda}^n)}
{d_0^{1,n-1}\left(\C\Diff(\kappa,\overline{\Lambda}^{n-1}(B,\NN))\right)}.
\end{equation}
This puts in evidence that the complex $E^1_0(B,\NN)$ has common
$n$--cocycles with $E_0^1$ but less $n$--coboundaries than it.\par
Take an element $\vartheta=[\square]_{\im d_{0,\textrm{rel}}}$ of
$\widehat{\kappa}(B,\NN)$, and, by using the $\C$--Green formula
for $\square$ (see \cite{CoAn}, 4.1.4), decompose $\vartheta$ as
the sum $\vartheta=h+\vartheta'$, being $h=[\square^*(1)]_{\im
d_{0,\textrm{rel}}}$ and $\vartheta'=[d\circ\square']_{\im
d_{0,\textrm{rel}}}$, with
$\square'\in\C\Diff(\kappa,\overline{\Lambda}^{n-1})$. Even though
$\square'$ is not uniquely determined by $\square$, so is
$d\circ\square'$. But $H(i)(\vartheta')$ is zero, so there exist
an unique $\theta'\in\widehat{\kappa}(\NN)$ such that
$\overline{\partial}(\theta')=\vartheta'$. This proves the
following
\begin{proposition}\label{BellaIdentificazione}
The map $\vartheta\mapsto\theta'$ splits the sequence
\REF{eqSeqEsCorE1bis}.
\end{proposition}
Thus we can identify $\vartheta$ with the pair $(\square^*(1),\
\theta')\in \widehat{\kappa}\oplus\widehat{\kappa}(\NN)$.

\subsection{The relative Euler operator}
\begin{definition}
The differential
$d_{1,\textrm{\normalfont{rel}}}^{0,n}:\overline{H}^n(B,\NN)\to
\widehat{\kappa}(B,\NN)$ is called the relative Euler operator and
is denoted by $\boldsymbol{E}_\textrm{\normalfont{rel}}$.
\end{definition}
By applying the relative Euler operator to a Lagrangian
$L=[\overline{\omega}]\in \overline{H}^n(B,\NN)$ one gets the pair
$(\ell_{\overline{\omega}}^*(1),\theta_{\overline{\omega}}')$
accordingly to Proposition \ref{BellaIdentificazione}. On the
other hand,  the extremality condition for the corresponding
variational problem is $\boldsymbol{E}_\textrm{rel}(L)=0$, i.e.,
$\ell_{\overline{\omega}}^*(1)=0$ and
$\theta_{\overline{\omega}}'=0$. The first of these two conditions
is the classical Euler--Lagrange equation corresponding to $L$,
while the second one we shall call the \emph{transversality
conditions} for the variational problem for $L$ with free boundary
(see \cite{CSpecSeqII}, 8.5) to be conform with the terminology in
the standard variational calculus.

\end{document}